\documentclass{amsart}

\newtheorem{theorem}{Theorem}[section]

\theoremstyle{definition}

\newtheorem{algorithm}[theorem]{Algorithm}
\newtheorem{conjecture}[theorem]{Conjecture}
\theoremstyle{remark}

\usepackage{amsmath,amssymb,amsthm,amscd}
\usepackage{amsbsy}
\usepackage{amsfonts}
\usepackage[mathscr]{eucal}
\usepackage[OT2,T1]{fontenc}
\usepackage[T3,T4]{fontenc}
\usepackage{latexsym}
\usepackage{enumerate}
\usepackage{graphicx}
\usepackage{color}	
\usepackage{epsfig}
\usepackage[numbers]{natbib}
\usepackage{listings}
\usepackage{psfrag}
\usepackage{verbatim}
\usepackage{url}
\usepackage{fancyhdr}

\newcommand{\Z}{\mathbb{Z}}
\newcommand{\Q}{\mathbb{Q}}
\newcommand{\C}{\mathcal{C}}

\newcommand{\Pt}[1]{\mathscr{#1}}
\DeclareSymbolFont{cyrls}{OT2}{wncyr}{m}{n}
\DeclareMathSymbol{\X}{\mathalpha}{cyrls}{"48}
\DeclareMathSymbol{\Y}{\mathalpha}{cyrls}{"55}
\DeclareMathSymbol{\Si}{\mathalpha}{cyrls}{"62}
\DeclareMathSymbol{\Sha}{\mathalpha}{cyrls}{"58}

\DeclareSymbolFont{cows}{T4}{fcr10}{m}{n}
\DeclareMathSymbol{\EF}{\mathalpha}{cows}{"84}

\DeclareSymbolFont{cheese}{T3}{tipa10}{m}{n}
\DeclareMathSymbol{\Em}{\mathalpha}{cheese}{"4D}
\DeclareMathSymbol{\Ef}{\mathalpha}{cheese}{"EA}
\DeclareMathSymbol{\Eg}{\mathalpha}{cheese}{"A5}

\begin{document}

\title[A binomial identity on the least prime factor of an integer]{A binomial identity on the least prime factor of an integer}

\author[S. Hambleton]{Samuel A. Hambleton}

\address{School of Mathematics and Physics, University of Queensland, St. Lucia, Queensland, Australia 4072}

\email{sah@maths.uq.edu.au}


\date{July 28, 2011.}


\begin{abstract}
An identity for binomial symbols modulo an odd positive integer $n$ relating to the least prime factor of $n$ is proved. The identity is discussed within the context of Pell conics.
\end{abstract}

\maketitle

\section{Introduction}

Many results exist on identities relating to binomial coefficients $\binom{m}{r}$ modulo $n$ where $n$ is an odd positive integer \cite{Dickson}. Granville \cite{Gran} has given new results concerning $\binom{m}{r} \pmod{p^q}$ where $p$ is prime, with a nice account of known results. Perhaps the most well known identity on factorials modulo $n$ is Wilson's theorem, which states that a positive integer $n$ is prime if and only if $(n-1)! \equiv -1 \pmod{n}$. Granville \cite{Gran} writes that Fleck \cite{Dickson} has generalized Wilson's theorem to the statement that for all positive integers $r$ less that the least prime divisor of $n$, $n$ is prime if and only if 
\begin{equation*}
\prod_{j=0}^{n-1-r} \binom{r+j}{r} \equiv (-1)^{\binom{r+1}{2}} \prod_{j=1}^{r- 1} \binom{r}{j} \pmod{n} .
\end{equation*}
Similarly, we will consider the residue modulo an odd positive integer $n$ of a symbol $\beta (n, r)$ defined in terms of binomial coefficients where, likewise, $r$ is less than or equal to the least prime divisor $p$ of $n$. We will briefly discuss the case $r > p$. Let $\lfloor a \rfloor$ and $\lceil a \rceil $ respectively denote the greatest integer $A \leq a$, and the least integer $A \geq a$.

\begin{theorem}\label{conics}
Let $n$ be an odd positive integer, let $r \geq 2$ be an integer, and let $p$ be the least prime divisor of $n$. Define $\alpha (n,r)$ to be the non-negative residue modulo $n$ of
\begin{equation}\label{thecoeff}
\beta (n, r) = (-1)^{\lfloor \frac{r}{2} \rfloor } \binom{\frac{n-1}{2}-\lceil \frac{r}{2} \rceil }{ \lfloor \frac{r}{2} \rfloor } - \binom{\frac{n-1}{2}}{r} (-2)^r .
\end{equation}
Then $\alpha (n,r)$ satisfies $\alpha(n,r) = \left\{\begin{array}{ll}  0 \pmod{n}   &\mbox{if $r < p$ } \\ n/p \pmod{n} &\mbox{if $r = p$ } \end{array}\right.$ .
\end{theorem}

Eqn. \eqref{thecoeff} occurs as the leading coefficient of the difference modulo $n$ of two polynomials which are important in the study of the affine genus zero curves known as Pell conics examined in detail by Lemmermeyer \cite{Lem03, MCP} and other authors \cite{LLPC, HambSchar} in relation to the analogy between these curves and elliptic curves. Let $\Delta$ be the fundamental discriminant of a quadratic number field $K = \Q(\sqrt{\Delta })$. Pell conics are the curves 
\begin{equation*}
\C : \X^2 - \Delta \Y^2 = 4 ,
\end{equation*}
with group law 
\begin{equation}\label{grplaw}
\Pt{P}_1 + \Pt{P}_2 = \Bigl( \frac{\X_1 \X_2 + \Delta \Y_1 \Y_2 }{2}, \frac{\X_1 \Y_2 + \X_2 \Y_1 }{2} \Bigr)
\end{equation}
defined for points $\Pt{P}_1 = (\X_1 , \Y_1)$ and $\Pt{P}_2 = (\X_2 , \Y_2)$ over $(\Z / n)$, $\Z$, $\Q$, and algebraic numbers $\overline{\Q}$ among various other rings $R$ for which the binary operation $+$ of Eqn. \eqref{grplaw} forms a group $\C(R)$ with identity $(2, 0)$. See \cite{Lem03} for more on these curves.

We define the polynomials $\EF_n({\X})$ by 
\begin{equation*}
\EF_1 = 1, \EF_3 = \X + 1, \EF_{2 j + 3} = \X \EF_{2 j +1} - \EF_{2 j -1},
\end{equation*} 
The origin of the polynomials $\EF_{n}(\X )$ can be traced to D. H. Lehmer \cite{Leh30} who has compared a Lucas function to Sylvester polynomials $\Psi_n(x, y)$ appearing in Bachmann's \cite{Sylvester} book. The polynomials $\Psi_n(x, y)$ correspond to the $G_m(x)$ used by Williams \cite{W87}.
\begin{equation*}
\EF_n(\X ) = G_{(n-1)/2}(\X ) \text{\em of Williams} = \Psi_{n}(\X , 1) \text{\em of Sylvester according to Lehmer}.
\end{equation*}

It has been shown \cite{LLPC, HambSchar} that the zeros of the polynomials $\EF_n(\X )$ are in one to one correspondence with the $\X$-coordinates of the non-trivial points $\Pt{P} \not= (2, 0)$ of order dividing $n$ in the group $\C (\overline{\Q})$, non-trivial points of the $n$-torsion subgroup $\C (\overline{\Q})[n]$. One simply expresses the $\X$-coordinate of $n(\X , \Y)$, meaning $n-1$ additions $(\X , \Y) + (\X , \Y) + \dots (\X , \Y)$, as $(\X - 2) \EF_n(\X)^2 + 2$. In order to give a proof of quadratic reciprocity \cite{HambSchar} using $p$-torsion on Pell conics where $p$ is an odd prime, it was demonstrated that
\begin{equation*}
\EF_p(\X) \equiv (\X - 2)^{\frac{p-1}{2}} \pmod{p} .
\end{equation*} 
The leading coefficient of the polynomial $\EF_n(\X) - (\X - 2)^{\frac{n-1}{2}}$ evaluated modulo $n$ is the more general question which we address. The polynomials $\EF_n$ are also discussed in the context of Dickson polynomials of the second kind, $E_n(x, a) = \sum_{j=0}^{\lfloor n / 2 \rfloor} \binom{n-j}{j}(-a)^j x^{n - 2 j}$. In particular, the identity, p.32 of \cite{LMT93}, 
\begin{equation*} 
\EF_{2 n + 1}(\X ) = E_{n}(\X , 1 ) + E_{n-1}(\X , 1 ) ,
\end{equation*} 
allows writing, for odd $n$, 
\begin{equation*} 
\EF_{n}(\X ) = \sum_{r=0}^{\frac{n-1}{2}} (-1)^{\lfloor \frac{r}{2} \rfloor } \binom{\frac{n-1}{2}-\lceil \frac{r}{2} \rceil }{ \lfloor \frac{r}{2} \rfloor }  {\X}^{\frac{n-1}{2}-r} .
\end{equation*}  
This completes the discussion of the context of the identity for $\beta (n, r)$.

\section{Proof of the main result}

We require the following equality which holds for all positive integers $a$.
\begin{equation}
\label{firstident} \prod_{j=1}^a (a+j) = 2^a \prod_{j = 0}^{a-1} (2 j +1) .
\end{equation}
Eqn. \eqref{firstident} may be proved by reordering the products in the numerator and denominator of $\prod_{j=1}^a \frac{a+j}{4 j -2}$, showing that this is equal to $1$. The proof of Theorem \ref{conics} is as follows.

\begin{proof}
First assume that $r < p$. Let $s = \lfloor r / 2 \rfloor$ and $t = \lceil r / 2 \rceil$. Then
\small
\begin{eqnarray*}
\beta (n,r) & = & (-1)^{s} \binom{\frac{n-1}{2} - t }{ s } - \binom{\frac{n-1}{2}}{r} (-2)^{r} , \\
           & = & \Bigl( \frac{(-1)^s }{s !} - \frac{(-2)^{r} \prod_{j=0}^{t -1} \bigl( \frac{n-1}{2} -j \bigr) }{r!} \Bigr) \prod_{j=1}^{t - 1} \bigl( \frac{n-1}{2} - s - j \bigr) , \\
           & = & \Bigl( \frac{(-1)^s \prod_{j=1}^{t} (s + j) - (-2)^{r} \prod_{j=0}^{t -1} \bigl( \frac{n-1}{2} -j \bigr) }{r!} \Bigr) \prod_{j=1}^{t -1} \bigl( \frac{n-1}{2} - s - j \bigr) , \\
           & = & \Bigl( \frac{(-1)^s \prod_{j=1}^{t} (s + j) - (-1)^{r+t} 2^s \prod_{j=0}^{t -1} \bigl( 1 + 2 j - n \bigr) }{r! (-2)^{t -1}} \Bigr) \prod_{j=1}^{t -1} \bigl( 1 + 2 s + 2 j - n \bigr) , \\
           & = & \Bigl( \frac{ \prod_{j=1}^{t} (s + j) - 2^s \prod_{j=0}^{t -1} \bigl( 1 + 2 j - n \bigr) }{r!} \Bigr) 2^{-t +1} (-1)^{r -1} \prod_{j=1}^{t -1} \bigl( 1 + 2 s + 2 j - n \bigr) , \\
\alpha (n, p) & \equiv & \Bigl( \frac{ \prod_{j=1}^{t} (s + j) - 2^s \prod_{j=0}^{t -1} \bigl( 1 + 2 j \bigr) }{r!} \Bigr) 2^{-t +1} (-1)^{r-1} \prod_{j=1}^{t -1} \bigl( 1 + 2 s + 2 j \bigr) \pmod{n} . 
\end{eqnarray*}
\normalsize
Since $r$ is strictly less than $p$, the integers $r!$ and $n$ are relatively prime. By Eqn. \eqref{firstident}, $\alpha (n, r) = 0$. Now let $r = p = 2 s + 1$. Then 
\small
\begin{eqnarray*}
\beta (n, p) & = & (-1)^{s} \binom{\frac{n-1}{2} - s - 1 }{ s } + \binom{\frac{n-1}{2}}{p} 2^p , \\
            & = & \Bigl( \frac{(-1)^s }{s !} +\frac{2^s \prod_{j=0}^s \bigl( n - 1 - 2 j \bigr) }{p!} \Bigr) \prod_{j=1}^s \bigl( \frac{n-1}{2} - s - j \bigr) , \\
            & = & \Bigl( \frac{(-1)^s }{s !} +\frac{2^s (\frac{n}{p} - 1) \prod_{j=0}^{s-1} \bigl( n - 1 - 2 j \bigr) }{(p-1)!} \Bigr) \prod_{j=1}^s \bigl( \frac{n-1}{2} - s - j \bigr) , \\
            & = & \frac{ \prod_{j=0}^{s-1} (s + j + 1) + 2^s \Bigl( \frac{n}{p} - 1 \Bigr) \prod_{j=0}^{s-1} \bigl( -n + 1 + 2 j \bigr) }{(p-1)! 2^{s}} \prod_{j=1}^s \bigl( -n + p + 2 j \bigr) , \\
\alpha (n, p) & \equiv & \Bigl( \prod_{j=0}^{s-1} (s + j + 1) + \Bigl( \frac{n}{p} - 1 \Bigr) \prod_{j=1}^{s} \bigl( s + j \bigr) \Bigr) (p-1)!^{-1} 2^{-s} \prod_{j=1}^s \bigl( p + 2 j \bigr) \pmod{n} , \\
            & \equiv & \frac{n}{p}  (p-1)!^{-1} 2^{-s} \prod_{j=1}^s \bigl( s + j \bigr) \bigl( p + 2 j \bigr) \pmod{n} , \\
            & \equiv & \frac{n}{p} (p-1)!^{-1} 2^{-p+1} \prod_{j=1}^{p-1} \bigl( p + j \bigr) \pmod{n} , \\
            & \equiv & \frac{n}{p} \prod_{j=1}^{p-1} (2 j)^{-1} \bigl( p + j \bigr) \pmod{n} .
\end{eqnarray*}
\normalsize
Fermat's theorem shows that $\prod_{j=1}^{p-1} (2 j)^{-1} \bigl( p + j \bigr) \equiv 1 \pmod{p}$. It follows that $\alpha (n , p) = \frac{n}{p}$.
\end{proof}

We conclude by speculating as to the value of $\alpha (n, r)$ when $r$ exceeds the least prime divisor of $n$, within some bounds. The author has only tested the following conjecture for $n < 10^6$.

\begin{conjecture}\label{conjj}
Let $p$ be the least prime divisor of an odd integer $n$ and assume that $2 \sqrt{n} < 3 p$. If $r$ is an integer bounded by $p < r < \sqrt{n}$ then $\alpha (n, r) > 0$.
\end{conjecture}

If Conjecture \ref{conjj} holds and the least prime divisor $p$ of $n$ satisfies $2 \sqrt{n} < 3 p$ then the follow exponential algorithm will terminate.

\begin{algorithm}
Let $A = (a_1, a_2)$ and assume we wish to factor $n$. Set $A = (2, \lfloor \sqrt{n} \rfloor ) $. If $\alpha \Bigl( n , \Bigl\lfloor \frac{a_1 + a_2}{2} \Bigr\rfloor \Bigr) = 0$, Set $A = \Bigl( \Bigl\lfloor \frac{a_1 + a_2}{2} \Bigr\rfloor , a_2\Bigr)$, otherwise set $A = \Bigl( a_1, \Bigl\lfloor \frac{a_1 + a_2}{2} \Bigr\rfloor \Bigr)$, and print $A$. Repeat until $a_2 - a_1 \leq 2$.
\end{algorithm}

\section*{Acknowledgments}

The author would like to thank Victor Scharaschkin for doctoral supervision of which this project has been a very small part of, and supported by the University of Queensland.

\bibliographystyle{amsplain}

\end{document}